\documentclass{swpaperproc}
\usepackage{latexsym}
\usepackage{graphicx}
\usepackage{mathptmx}

\usepackage{amsmath}
\usepackage{amsfonts}
\usepackage{amssymb}
\usepackage{amsbsy}
\usepackage{amsthm}

\usepackage[pdftex,colorlinks=true,urlcolor=blue,citecolor=black,anchorcolor=black,linkcolor=black]{hyperref}

\usepackage{dirtytalk}

\newtheoremstyle{sw}
{3pt}
{3pt}
{}
{}
{\bf}
{}
{.5em}
{}

\theoremstyle{sw}

\begin{document}


%
%

\pagestyle{fancyplain}

\thispagestyle{plain}
\firstPageHead{}

\chead{\fancyplain{}{\itshape Amini and Van Nieuwenhuyse}}

\rhead{}
\cfoot{}
\renewcommand{\headrulewidth}{0pt} 

\makeatletter
\let\@internalcite\cite
\def\cite{\def\@citeseppen{-1000}%
    \def\@cite##1##2{(##1\if@tempswa , ##2\fi)}%
    \def\citeauthoryear##1##2##3{##1 ##3}\@internalcite}
\def\citeNP{\def\@citeseppen{-1000}%
    \def\@cite##1##2{##1\if@tempswa , ##2\fi}%
    \def\citeauthoryear##1##2##3{##1 ##3}\@internalcite}
\def\citeN{\def\@citeseppen{-1000}%
    \def\@cite##1##2{##1\if@tempswa, ##2)\else{}\fi}%
    \def\citeauthoryear##1##2##3{##1 (##3)}\@citedata}
\def\citeA{\def\@citeseppen{-1000}%
    \def\@cite##1##2{(##1\if@tempswa , ##2\fi)}%
    \def\citeauthoryear##1##2##3{##1}\@internalcite}
\def\citeANP{\def\@citeseppen{-1000}%
    \def\@cite##1##2{##1\if@tempswa , ##2\fi}%
    \def\citeauthoryear##1##2##3{##1}\@internalcite}
\def\shortcite{\def\@citeseppen{-1000}%
    \def\@cite##1##2{(##1\if@tempswa , ##2\fi)}%
    \def\citeauthoryear##1##2##3{##2 ##3}\@internalcite}
\def\shortciteNP{\def\@citeseppen{-1000}%
    \def\@cite##1##2{##1\if@tempswa , ##2\fi}%
    \def\citeauthoryear##1##2##3{##2 ##3}\@internalcite}
\def\shortciteN{\def\@citeseppen{-1000}%
    \def\@cite##1##2{##1\if@tempswa, ##2\else{}\fi}%
    \def\citeauthoryear##1##2##3{##2 (##3)}\@citedata}
\def\shortciteA{\def\@citeseppen{-1000}%
    \def\@cite##1##2{(##1\if@tempswa , ##2\fi)}%
    \def\citeauthoryear##1##2##3{##2}\@internalcite}
\def\shortciteANP{\def\@citeseppen{-1000}%
    \def\@cite##1##2{##1\if@tempswa , ##2\fi}%
    \def\citeauthoryear##1##2##3{##2}\@internalcite}
\def\citeyear{\def\@citeseppen{-1000}%
    \def\@cite##1##2{(##1\if@tempswa , ##2\fi)}%
    \def\citeauthoryear##1##2##3{##3}\@citedata}
\def\citeyearNP{\def\@citeseppen{-1000}%
    \def\@cite##1##2{##1\if@tempswa , ##2\fi}%
    \def\citeauthoryear##1##2##3{##3}\@citedata}
%
%
%
\def\@citedata{%
    \@ifnextchar [{\@tempswatrue\@citedatax}%
                  {\@tempswafalse\@citedatax[]}%
}

\def\@citedatax[#1]#2{%
\if@filesw\immediate\write\@auxout{\string\citation{#2}}\fi%
  \def\@citea{}\@cite{\@for\@citeb:=#2\do%
    {\@citea\def\@citea{, }\@ifundefined
       {b@\@citeb}{{\bf ?}%
       \@warning{Citation `\@citeb' on page \thepage \space undefined}}%
{\csname b@\@citeb\endcsname}}}{#1}}%

%
\def\@citex[#1]#2{%
\if@filesw\immediate\write\@auxout{\string\citation{#2}}\fi%
  \def\@citea{}\@cite{\@for\@citeb:=#2\do%
    {\@citea\def\@citea{, }\@ifundefined
       {b@\@citeb}{{\bf ?}%
       \@warning{Citation `\@citeb' on page \thepage \space undefined}}%
{\csname b@\@citeb\endcsname}}}{#1}}%

%
\def\@biblabel#1{}
\makeatother



\newdimen\bibindent
\bibindent=0.0em
\def\thebibliography#1{\section*{\refname}\list
   {}{\settowidth\labelwidth{[#1]}
   \leftmargin\parindent
   \itemindent -\parindent
   \listparindent \itemindent
   \itemsep 0pt
   \parsep 0pt}
   \def\newblock{}
   \sloppy
   \sfcode`\.=1000\relax}

\setlength{\baselineskip}{12.7pt}

\title{A TUTORIAL ON KRIGING-BASED STOCHASTIC SIMULATION OPTIMIZATION}

\author{\textit{Sasan Amini}\\ [11pt]
Hasselt University\\
Agoralaan, 3590 Diepenbeek, Belgium\\
sasan.amini@uhasselt.be\\
\and
\textit{Inneke Van Nieuwenhuyse}\\[11pt]
 Hasselt University\\
Agoralaan, 3590 Diepenbeek, Belgium\\
inneke.vannieuwenhuyse@uhasselt.be\\
 }

\maketitle

\section*{ABSTRACT}
This tutorial focuses on kriging-based simulation optimization, emphasizing the importance of data efficiency in optimization problems involving expensive simulation models. It discusses how kriging models contribute to developing algorithms that minimize the number of required simulations, particularly in the presence of noisy evaluations. The tutorial compares the performance of kriging-based algorithms against traditional polynomial-based optimization methods using an illustrative example. Additionally, it discusses key extensions of kriging-based algorithms, including multi-objective and constrained optimization, providing insights into their application in complex, real-world settings.

\vspace{3mm}
\noindent \textbf{Keywords:} Simulation Optimization, Metamodel, Kriging, Data Efficient Global Optimization

\section{INTRODUCTION}
\label{sec:intro}

Simulation optimization refers to the process of finding optimal decisions in complex systems where analytical methods are not applicable, and outcomes can only be evaluated through (often stochastic) simulation models \shortcite{fu2015handbook,amaran2016simulation,nelson2021foundations}. Classical stochastic simulation optimization methods, such as the genetic algorithms used in \shortciteN{ding2005simulation} and \shortciteN{pasandideh2006multi}, are computationally intensive as they require the evaluation of a large number of potential solutions. Additionally, multiple replications at each solution are needed to improve the accuracy of performance estimates, further increasing the computational burden.  In this tutorial, our focus is on problems for which the simulation model is expensive to evaluate, such that the analyst can only afford a limited number of evaluations. Furthermore, the simulation model is a black box where only input/output (I/O) information is available (we have no access to, for instance, derivative information that might guide the optimization).

To solve such challenging problems, researchers have increasingly relied on metamodels \shortcite{barton2023tutorial} as a way to approximate the behavior of simulation models. Metamodels, also known as surrogate models \shortcite{gramacy2020surrogates,hong2021surrogate}, act as simplified representations of the underlying simulation model, enabling efficient exploration of the solution space without the need for excessive computational resources. The use of metamodels has been widely studied in the literature, particularly for problems where simulation costs are high, making them invaluable in industrial settings, see, e.g., \shortciteN{stavropoulos2023metamodelling},  \shortciteN{pietrusewicz2019metamodelling}, and \shortciteN{hassannayebi2019train}. For a systematic literature review of metamodeling-based simulation optimization, we refer to  \shortciteN{hong2021surrogate} and \shortciteN{do2022metamodel}.

Traditional metamodeling approaches often rely on linear or quadratic regression models \shortcite{kleijnen2018design}. While these models are computationally inexpensive and simple to implement, they often lack the flexibility required to accurately capture the complex, nonlinear I/O relationships that are common in many real-world systems \shortcite{barton2023metamodelling}. This reduces their effectiveness, particularly when the simulation response surface exhibits significant curvature or interactions among variables that cannot be captured by simple polynomial terms.
In more recent years, \emph{kriging} metamodels have emerged as a powerful and versatile tool in simulation optimization \shortcite{do2024adaptive}. Kriging models, also known as Gaussian process regression (GPR) models in the field of machine learning and statistics \shortcite{shahriari2015taking,gramacy2020surrogates}, do not only provide accurate predictions of system performance but also estimate the uncertainty related to these predictions, referred to as the kriging variance. This is a key advantage, which has proven useful in \emph{data-efficient global optimization}, particularly since the seminal work of \shortciteN{jones1998efficient}. 
In engineering and machine learning, this area of research is often referred to as Bayesian optimization \shortcite{garnett2023bayesian}.

\section{Kriging models}

This method, originally developed for approximating I/O relationships in mining geology and named after South African engineer Danie G. Krige \shortcite{krige1951statistical}, has later been extended to computer simulation experiments by \shortciteN{10.1214/ss/1177012413}. Kriging is a spatial interpolation method used to obtain predictions at unsampled locations based on observed data \shortcite{kleijnen2009kriging}. It assumes that the distance between sample points reflects a spatial correlation that can be used to explain variation in the surface, i.e., the underlying unknown function \shortcite{chiles2018fifty}.

\subsection{Kriging for deterministic responses}
Kriging models can approximate the input-output (I/O) relationships of deterministic simulation models, i.e. models for which it can be assumed that the evaluation of the system for a given set of input values results in an \emph{exact} (i.e., non-noisy) observation of the corresponding output. In other words, if you were to repeat the evaluation at the same point, you would always get the same result; there is no randomness or variability on the output for a given input.

Ordinary kriging models the output of the unknown response function at a given input vector $\mathbf{x}$ as the sum of a constant term, and an observation from a zero-mean covariance-stationary Gaussian random field:

\begin{equation}
  f(\mathbf{x})=\beta_0+ M(\mathbf{x})
 \label{equ:0}
\end{equation}

\noindent where $\beta_0$ then represents the overall surface mean. In its most general form, this first term may also consist of a combination of prespecified functions at $\mathbf{x}$; this is referred to as \emph{universal kriging} \shortcite{palar2017multi}. Yet, the modeling assumption with the constant term is most common in the literature. 

Given the modeling assumption, the resulting kriging model yields a probability distribution over functions, such that for any finite set of input points, the corresponding function values are jointly Gaussian distributed \shortcite{gramacy2020surrogates,rasmussen2006gaussian}.
The Gaussian random field is uniquely characterized by its mean $m(.)$, which is assumed to be zero, and its covariance function $k(.)$, also referred to as \emph{kernel}. This covariance function models the pairwise covariances between function values (outputs) at different input locations. The basic assumption here is that the output values of two distinct points in the solution space are spatially correlated \shortcite{rasmussen2006gaussian}; i.e., the correlation between the output values depends on the distance between the input locations. Data points that are close to each other (in the solution space) are more likely to have similar output values, so they tend to have high correlation. As the spatial separation between data points increases, the correlation decreases. Many different types of covariance functions can be used; the choice of covariance function directly impacts the model's ability to capture the underlying function \shortcite{archetti2019bayesian,garnett2023bayesian,wang2023recent}.  
 Commonly used kernel functions in the literature are the squared exponential (Equation \ref{kernel se}), Matérn$_{3/2}$ (Equation \ref{kernel mat3/2}), and Matérn$_{5/2}$ (Equation \ref{kernel mat5/2}). 
 
\begin{figure}[t!]
\centering
\includegraphics[width=0.85\textwidth]{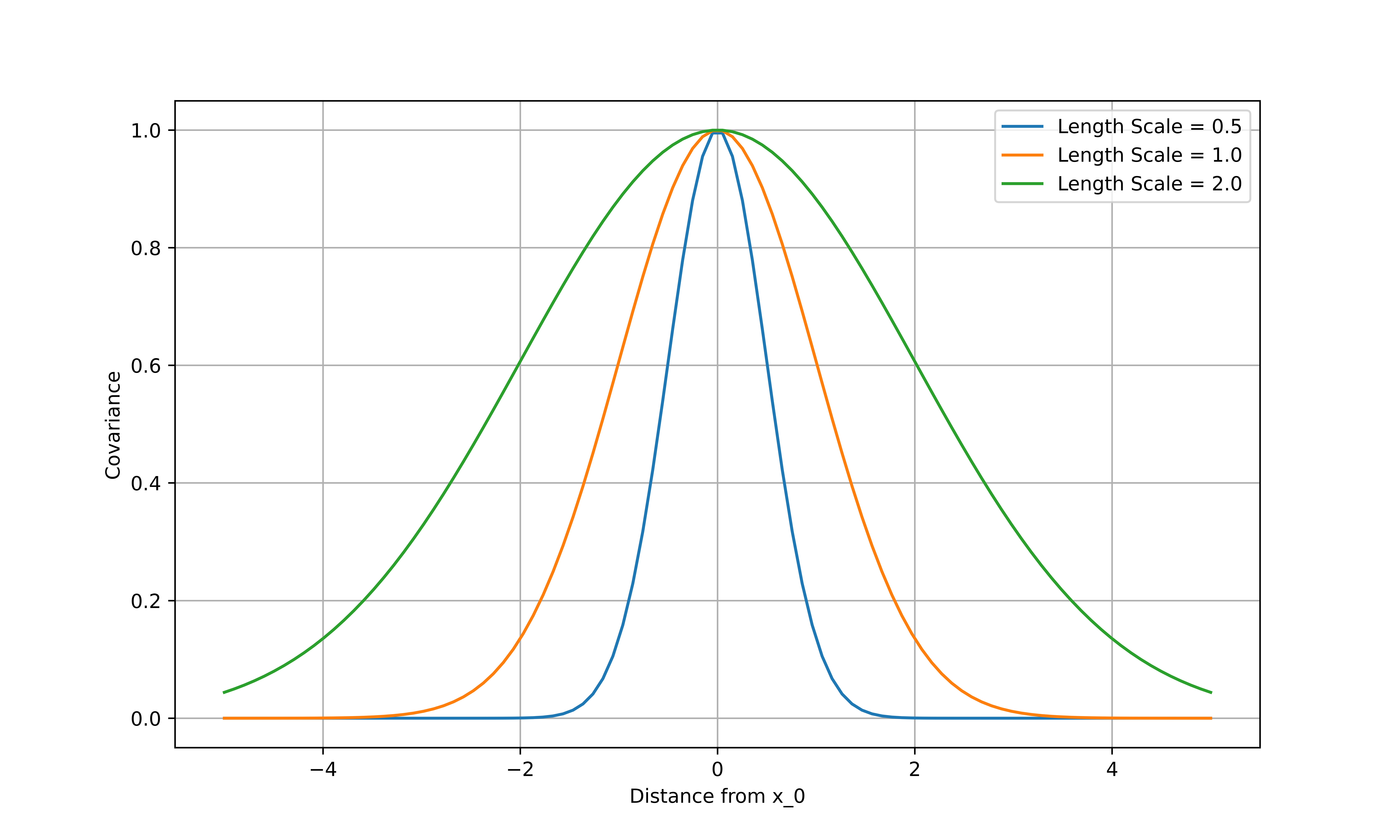}
\caption{Effect of distance and length-scale on spatial correlation using a squared exponential covariance function, assuming $\sigma^2=1$}
\label{fig: spatial correl}
\end{figure}

 In Equations \eqref{kernel se}, \eqref{kernel mat3/2}, and \eqref{kernel mat5/2}, $d$ is the Euclidean distance between the two input locations, $\sigma^2$ is the variance of the kriging model, and $l$ is a hyperparameter known as the \emph{length-scale}, which controls how quickly the correlation between the output of two input points decreases as they move further away in the solution space. 
 Figure \eqref{fig: spatial correl} illustrates the concept of spatial correlation using the squared exponential kernel. The x-axis represents the distance from a reference point $x_0$, while the y-axis shows the resulting covariance between the output values (for ease of interpretation, the plot assumes $\sigma^2 = 1$; hence, the covariance values are equal to the correlation values). Clearly, a smaller length-scale implies that the correlation drops more rapidly with distance, which may cause function values to change quickly over short distances. By contrast, a larger length-scale implies that the correlation decreases more gradually, implying smoother variations in the data \shortcite{rasmussen2006gaussian}.

\begin{equation}
    k_{SE}(.,.)=\sigma^2exp(-\frac{d^2}{2l^2})
    \label{kernel se}
\end{equation}

\begin{equation}
  k_{M_{3/2}}(.,.)=\sigma^2\bigl(1+\frac{\sqrt{3}d}{l}\bigr)exp\bigl(-\frac{\sqrt{3}d}{l}\bigr)
    \label{kernel mat3/2}  
\end{equation}

\begin{equation}
  k_{M_{5/2}}(.,.)=\sigma^2\bigl(1+\frac{\sqrt{5}d}{l}+\frac{5d^2}{3l^2}\bigr)exp\bigl(-\frac{\sqrt{5}d}{l}\bigr)
    \label{kernel mat5/2}  
\end{equation}

Ideally, the properties of the chosen covariance function, such as smoothness assumptions and the ability to capture specific patterns, should align with the characteristics of the function under study \shortcite{archetti2019bayesian}. This, however, requires an understanding of the problem domain and prior knowledge about the underlying function. Often, this information is not available. In such cases, kernel choice relies on expert knowledge and involves hyperparameter tuning techniques such as maximum likelihood estimation or cross-validation \shortcite{duvenaud2014kernel,duvenaud2013structure}. Some authors \shortcite{roman2014dynamic,roman2019experimental} also study the dynamic change of kernel functions during the optimization process using their proposed kernel selection
criteria.

For a chosen covariance function $k(.)$ and a set of $p$ observed points, $\mathbf{x}_1, ..., \mathbf{x}_p$, we can calculate the $p \times p$ covariance matrix:

\begin{equation}
K=
\left (
\begin{array}{ccc}
\begin{array}{l}
k(\mathbf{x}_1,\mathbf{x}_1)\\
\end{array}
& \cdots & 
\begin{array}{l}
k(\mathbf{x}_1,\mathbf{x}_p) \\
\end{array} \\
\vdots & \ddots & \vdots\\
\begin{array}{l}
k(\mathbf{x}_p,\mathbf{x}_1) \\
\end{array} &
\cdots & 
\begin{array}{l}
k(\mathbf{x}_p,\mathbf{x}_p) \\
\end{array} 
\end{array}
\right )
\end{equation}

\noindent The predictor value of the unknown function $f(\mathbf{x})$ at an arbitrary point $\mathbf{x}_{p+1}$ is given by

\begin{equation}\label{mu}
    \mu(\mathbf{x}_{p+1})=\hat{\beta_0}+\mathbf{k}^TK^{-1}(f_{1:p}-\mathbf{1}\hat{\beta_0})
\end{equation}

\noindent and the variance (also referred to as mean squared error, MSE) of this predictor by

\begin{equation}\label{sigma}
    s^2(\mathbf{x}_{p+1})= k(\mathbf{x}_{p+1}, \mathbf{x}_{p+1}) -\mathbf{k}^TK^{-1}\mathbf{k}+\frac{(1-\mathbf{1}^TK^{-1}\mathbf{k})^2}{\mathbf{1}^TK^{-1}\mathbf{1}}
\end{equation}

\noindent where $\mathbf{k}=\bigr[k(\mathbf{x}_{p+1},\mathbf{x}_1), k(\mathbf{x}_{p+1},\mathbf{x}_2), \dots, k(\mathbf{x}_{p+1},\mathbf{x}_p)\bigl]$, and $f_{1:p}$ is the matrix of the observed output values so far. In Equations \eqref{mu} and \eqref{sigma}, $\hat{\beta_0}$ is an estimation of $\beta_0$, the overall surface mean and can be calculated by:

\begin{equation}
\hat{\beta_0}=\frac{\mathbf{1}^Tk^{-1}f_{1:p}}{\mathbf{1}^TK^{-1}\mathbf{1}}
\end{equation}

The last term in Equation \eqref{sigma} which adjusts the uncertainty to account for the estimation of the constant mean, $\beta_0$, drops when $\beta_0$ will be considered as a given parameter and known.

The key modeling assumption implies that the joint distribution of the observed outputs and the unknown output at a new input location ($\mathbf{x}_{p+1}$) is a multivariate normal distribution. The predictive distribution of the unknown output at the individual location $\mathbf{x}_{p+1}$ is thus normally distributed, with mean equal to the predictor value, and variance given by the MSE. As will be shown in section \ref{opt}, this distribution information significantly aids the exploration-exploitation balance within the context of efficient global optimization \shortcite{garnett2023bayesian}. As the output observations are treated as \emph{deterministic}, the kriging model interpolates between these observations: at simulated input/output locations, the predictor value coincides with the simulated outputs, and the MSE of the predictor is zero. 

\begin{figure}[t!]
\centering
\includegraphics[width=\textwidth]{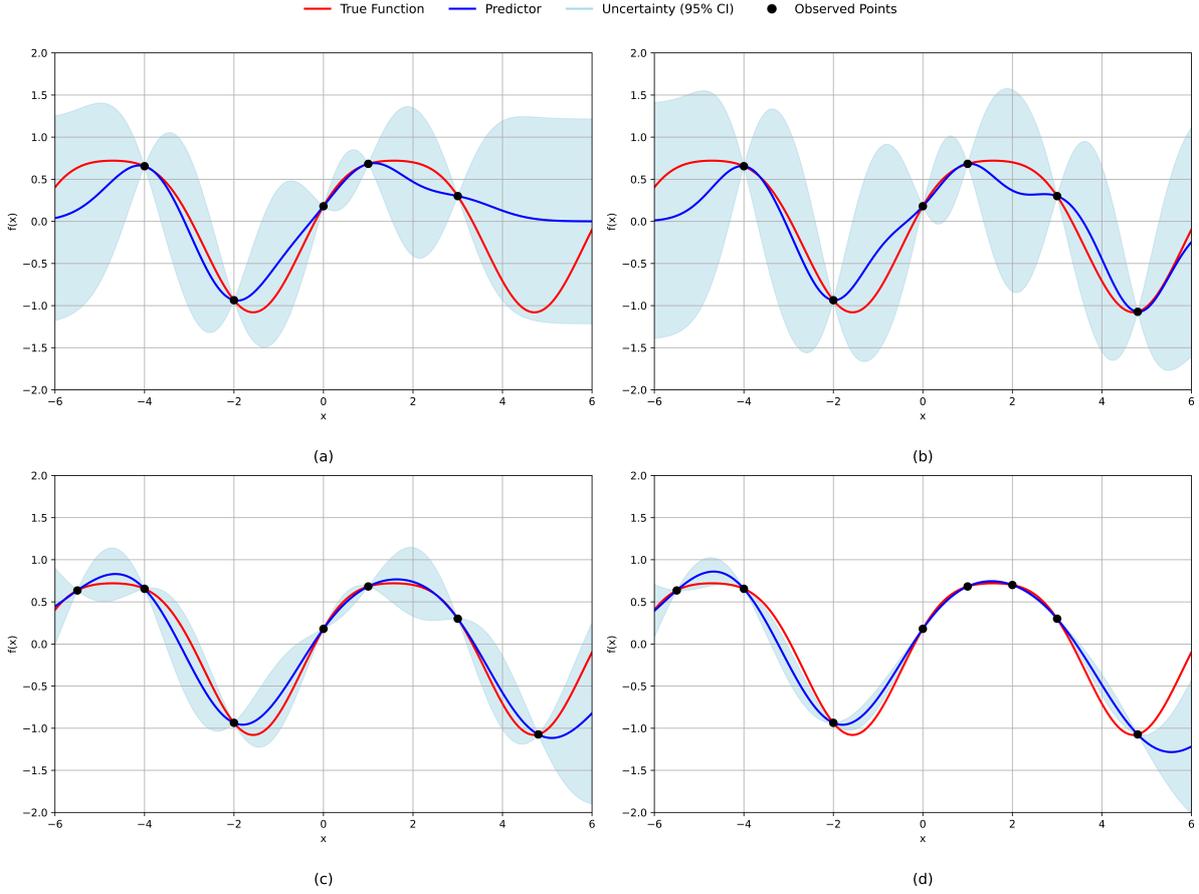}
\caption{Progressive  improvement of kriging model predictions as additional data points are added}
\label{fig: kriging improv}
\end{figure}

Figure \eqref{fig: kriging improv} demonstrates how the kriging model progressively improves its predictions as additional data points are incorporated. In the top-left plot, the model starts with only five observed points, showing relatively high uncertainty in regions without data. In the top-right plot, an additional point is added, reducing uncertainty around the corresponding region and improving prediction accuracy. The plots in the bottom row show how observing more points further refines the model's predictions. 

Note also that the model’s uncertainty (represented by the blue-shaded regions) drops to zero at the observed points, as the observed outputs are treated as deterministic. 

\subsection{Kriging for stochastic responses: stochastic kriging}

\begin{figure}[t!]
\centering
\includegraphics[width=\textwidth]{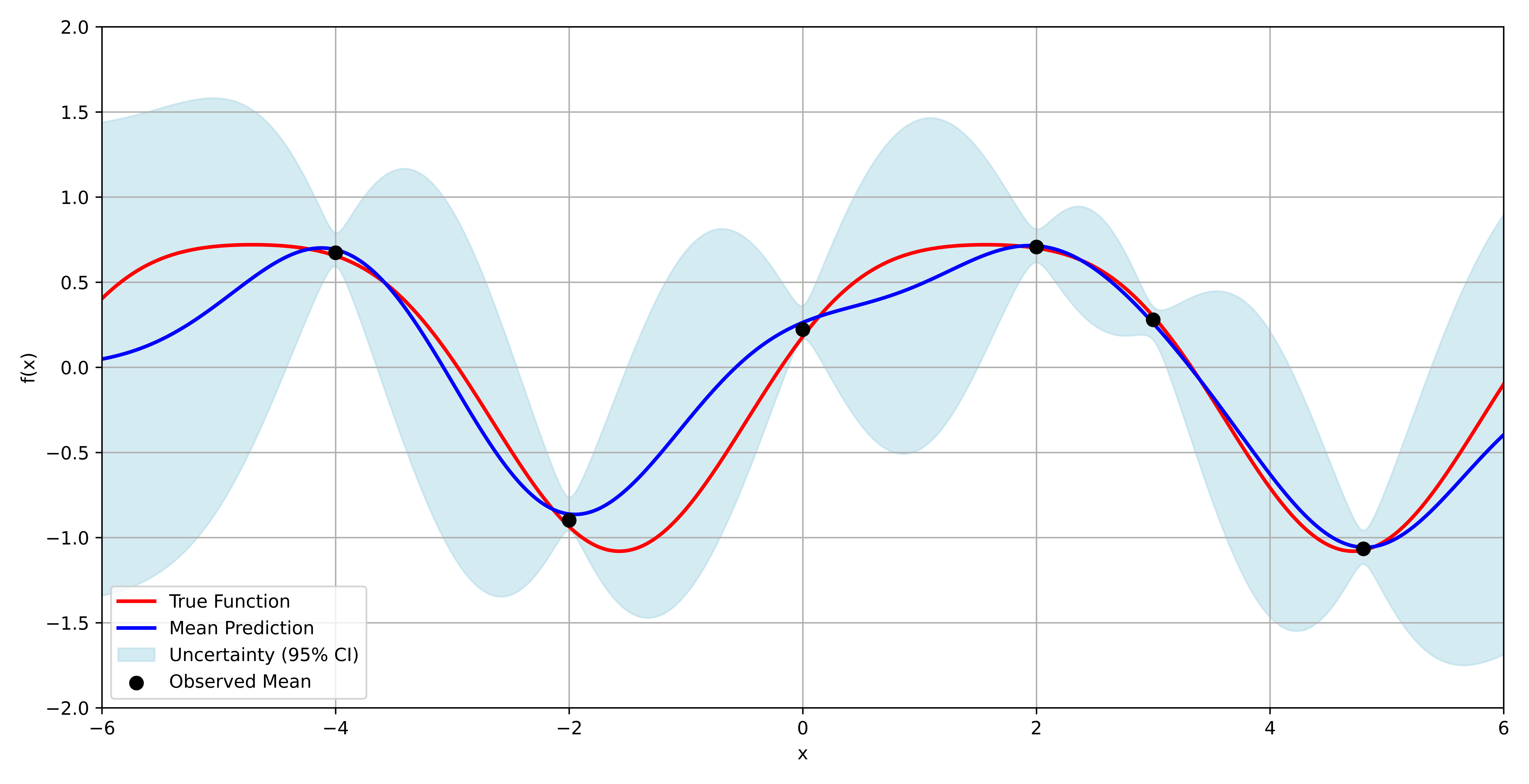}
\caption{Non-zero uncertainty at observed points in stochastic kriging}
\label{fig: sk}
\end{figure}

In real-world simulations, we often deal with stochastic systems: evaluating the system at the same point multiple times yields different results, as the system’s behavior is influenced by inherent randomness. Consequently, the output observations are noisy, and multiple independent evaluations (referred to as \emph{replications}) are typically required to estimate system performance reliably. In real-life systems, this noise is often \emph{heteroscedastic}, meaning that it differs depending on the input location. 

One way to account for this uncertainty in the outputs is the use of Stochastic Kriging \shortcite{ankenman_nelson_staum_2010}.  For other methods and more details, we refer to \shortciteN{binois2018practical}. 
Stochastic kriging assumes that, at an arbitrary design point $\mathbf{x}$, the function value obtained in the $\emph{r}^{th}$ replication can be modelled as:

\begin{equation}
 {f}_r(\mathbf{x})= \beta_0 + M(\mathbf{x}) + \epsilon_r (\mathbf{x})
 \label{equ:2}
\end{equation}

\noindent The first two terms of this equation are the same as for deterministic responses. The term $\epsilon_r (\mathbf{x})$ is referred to as \emph{intrinsic uncertainty}, as it is the uncertainty inherent in stochastic simulation. The intrinsic uncertainty is (naturally) independent and identically distributed across replications, having mean 0 and variance
$\tau^2(\mathbf{x})$ at any arbitrary point $x_i$. Note that the model allows for heteroscedastic noise, implying $\tau^2(\mathbf{x})$ need \emph{not} to be constant throughout the design space.
The stochastic kriging prediction at an arbitrary location $\mathbf{x}$ (whether observed or not) is given by:

\begin{equation}
 \mu_{sk}(\mathbf{x})= \hat{\beta_0}+\mathbf{k}^T[K+K_\epsilon]^{-1}(\overline{f}_{1:p}-\mathbf{1}\hat{\beta_0})
 \label{equ:3}
\end{equation}

\noindent where $\overline{f}$ is the $p\times 1$ vector containing the sample means obtained at the $p$ observed design points. $K_\epsilon$ is a diagonal matrix of size $p\times p$, showing the variance of the sample means on the main diagonal. $K_\epsilon$ thus reflects the intrinsic uncertainty of the system, as it is caused by the (heteroscedastic) noise in the output observations.
The mean squared error (MSE) of this predictor is given by:

\begin{equation}
 s_{sk}^{2}(\mathbf{x})= k(\mathbf{x},\mathbf{x}) - \mathbf{k}^T[K+K_\epsilon]^{-1}\mathbf{k}+ \frac{\Gamma^T\Gamma}{{1_p}^T[K+K_\epsilon]^{-1}{1_p}}
 \label{equ:4}
\end{equation}

\noindent where $\Gamma=1- {1_p}^T[K+K_\epsilon]^{-1}\mathbf{k}$. The essential difference between stochastic kriging and ordinary kriging is the presence of $K_\epsilon$ in the expressions for the predictor and the MSE. In the absence of noise, $K_\epsilon$ disappears from the equations, and the stochastic kriging predictor and MSE expressions boil down to the expressions of the ordinary kriging model.

Figure \eqref{fig: sk} illustrates the behavior of stochastic kriging in a setting with noisy observations. Contrary to ordinary kriging, which interpolates between observed points, the stochastic kriging predictor does not coincide with the sample means at the observed locations, and the corresponding MSE reflects the intrinsic uncertainty.


\section{Data-efficient simulation optimization with kriging models}\label{opt}

\begin{figure}[t!]
\centering
\includegraphics[width=0.85\textwidth]{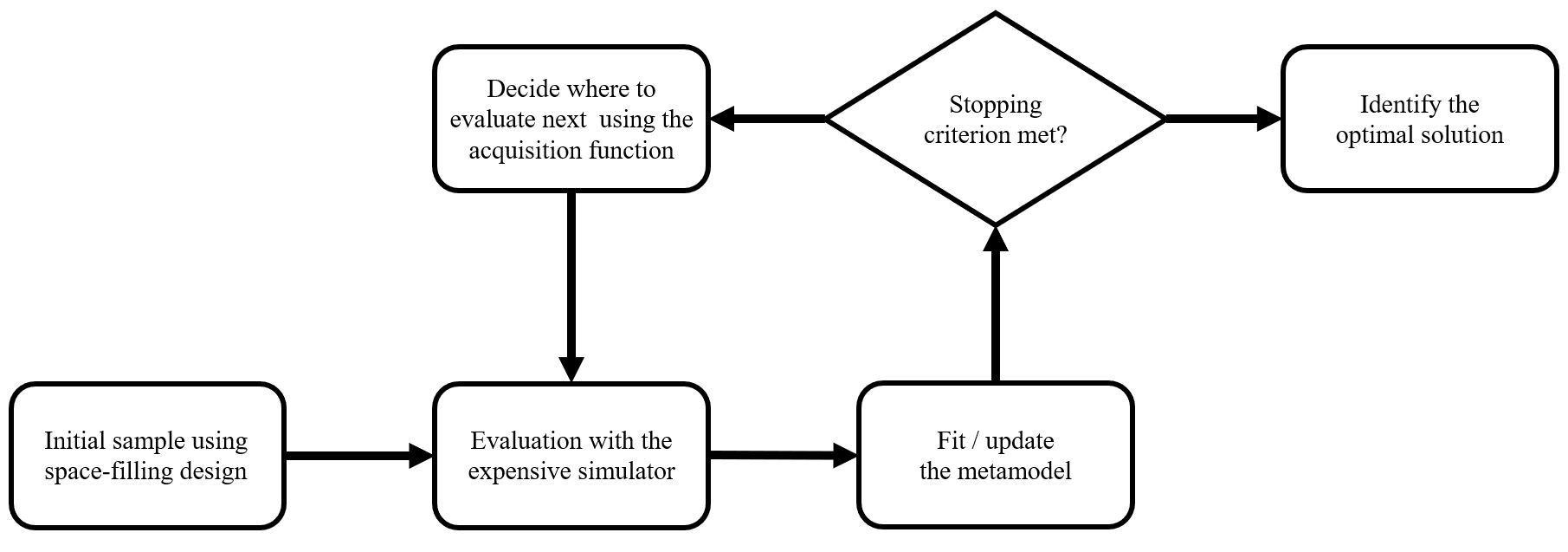}
\caption{General steps in a kriging-based simulation optimization algorithm}
\label{fig: BO steps}
\end{figure}

Kriging-based simulation optimization is a sequential optimization method that is particularly well suited for problems involving expensive, potentially nonconvex black box functions and a limited computational budget \shortcite{garnett2023bayesian,wang2023recent,frazier2018bayesian,candelieri2021gentle}. 
Various algorithms have been developed, yet they all share the same key steps, shown in Figure \ref{fig: BO steps}. In what follows, we discuss these steps in more detail.

\begin{itemize}

\item{\textbf{Initial sample:}} The first step is the generation of an initial sample. It has been shown that the success of metamodel-based optimization algorithms in achieving fast convergence to the optimal solution largely depends on the quality of the initial data used to train the metamodels \shortcite{vu2017surrogate}. The purpose of this step is to create a space-filling set of design points, i.e. a set of points that effectively covers the entire space without leaving significant gaps or underrepresented regions \shortcite{greenhill2020bayesian}. 
Various mathematical techniques and algorithms, such as Latin hypercube sampling \shortcite{jones1998efficient} and Sobol sequences \shortcite{picheny2013benchmark}, have been employed in the literature to generate space-filling designs. 

When repeatedly solving similar optimization problems, the analyst may learn how to set the initial design for the optimization algorithm from experience. Some authors have proposed sophisticated \emph{warm start} initialization procedures for hyperparameter tuning; see \shortcite{feurer2015initializing} for more details.
In general, the Latin Hypercube Sampling (LHS) design \shortcite{park1994optimal} is the most often used method in the Bayesian Optimization (BO) literature, given its space-filling and non-collapsing features. For more details on initial sampling techniques for surrogate models, we refer to \shortciteN{vu2017surrogate}.

\item{\textbf{Evaluation:}} Information about the unknown function(s) is acquired by running the expensive simulator at the selected points.
In stochastic simulations, variability is inherent in the simulation\shortcite{nelson2021foundations}. Increasing the number of replications improves the accuracy of the sample means in such settings, yet, as the total number of evaluations is limited, there is a trade-off between the number of simulations spent to replicating already observed points, and those spent to evaluate new points in the design space.

\item{\textbf{Build/Update Metamodel:}} 
In this step, the kriging metamodel is fit to the available I/O observations, such that the belief about the unknown function is updated \shortcite{frazier2018bayesian}. 


\item{\textbf{Search:}}  
Kriging-based simulation optimization is a sequential process, so after evaluating a new design point, the collected data is incorporated into the existing data, and the metamodel is updated. 
In this step, an \emph{acquisition function} (also referred to as \emph{infill criterion} \shortcite{rojas2020survey}) is used to determine the next solution to evaluate. This criterion quantifies the desirability of simulating a new solution \shortcite{frazier2018bayesian,archetti2019bayesian} based on the metamodel information. 

The uncertainty measure is crucial for guiding the search for optimal solutions, helping the optimization algorithm identify areas of the solution space that require further exploration. By combining the kriging predictions and the kriging variance, the acquisition function enhances the efficiency of simulation optimization, particularly in data-scarce environments where simulation runs are costly \shortcite{kleijnen2018design}. It typically trades off exploitation of the search space (i.e., sampling in areas with promising predictor values) versus exploration (sampling in areas with high predictor uncertainty
). 
During the early stages of optimization, exploration is often prioritized to gather diverse samples and gain a better understanding of the search space. This allows the algorithm to build an initial surrogate model that estimates the overall behavior of the function. As the optimization progresses, the surrogate model becomes more accurate, and the algorithm shifts towards exploitation \shortcite{gupta2006interplay,gan2021acquisition}. This allows for a more targeted search, potentially leading to a more accurate surrogate model in the neighborhood of the optimal solution and, as a result, faster convergence.
%
This iterative process of model updating and acquisition function maximization continues until the stopping criterion is met. \\
The choice of acquisition function depends on the specific problem characteristics, and is further discussed in Section \ref{sec:AF}. While some of these may outperform others in specific settings, to the best of our knowledge, there is no evidence that any of them consistently outperform others. 

\item{\textbf{Termination:}} At some moment, the optimization algorithm has to stop. Usually, this happens when it has reached a satisfactory solution \shortcite{Nguyen2023bayesian},  when it has depleted the available budget \shortcite{garnett2023bayesian}, or when the expected gain from the optimal observation is no longer worth the cost of the evaluation \shortcite{archetti2019bayesian}.

\item{\textbf{Identification:}} 
When the stopping criterion is met, the algorithm needs to identify the optimal solution(s). 
In the current literature, the main focus is on finding a solution that maximizes the expected performance of the objective function of interest, and it is common to identify the solution with the best-observed performance (i.e., the best function outcome in case of deterministic simulation, and the best sample mean in case of stochastic simulation) as the optimum \shortcite{carpio2018enhanced}.

\end{itemize}

\section{Acquisition functions}\label{sec:AF}

The acquisition function, combined with the metamodel, forms the \say{heart} of every kriging-based optimization algorithm \shortcite{garnett2023bayesian}. As mentioned above, it determines the next point to evaluate, given the current probabilistic belief about the function behaviour (i.e., the kriging model). 
In deterministic settings, Probability of Improvement (PI) and Expected Improvement (EI) are undoubtedly among the most popular functions. PI \shortcite{kushner1964new} is evaluated as:

\begin{equation}\label{Pr of Im}
    PI(\mathbf{x})=Pr\bigl(f(\mathbf{x}) \leq f_{min}) \bigr)=\Phi\Bigl(\frac{f_{min}-\mu(\mathbf{x})}{s(\mathbf{x})} \Bigr)
\end{equation}

\noindent where $f_{min}$ represents the best  objective function value observed so far, and $\Phi$ denotes the normal cumulative distribution.
Expected Improvement (EI), first proposed by \shortciteN{jones1998efficient}, has been widely used in the literature and has triggered the development of many variants. EI estimates the \emph{magnitude} of the improvement over the current best observed value, and is given by:  

\begin{equation}
\begin{aligned}
{\textstyle
EI(\mathbf{x})= \left[f_{min}- \mu(\mathbf{x})\right]\Phi \left(\frac{f_{min}- \mu(\mathbf{x})}{s(\mathbf{x})}\right)+ s(\mathbf{x}) \mathbf{\phi}\left(\frac{f_{min}- \mu(\mathbf{x})}{s(\mathbf{x})}\right)}
\end{aligned}
\label{equ:ei}
\end{equation}

\noindent where $\mathbf{\phi}$ denotes the normal probability density function.
Other popular acquisition functions in deterministic settings are confidence bounds, where upper and lower are used, respectively, for maximization
and minimization problems.

The acquisition functions discussed so far are myopic, in the sense that they focus on immediate improvements without considering the long-term impact of a sampling decision. In contrast, the Knowledge Gradient (KG) acquisition function evaluates the expected value of information that a new sample will provide \shortcite{garnett2023bayesian}. It focuses on maximizing the expected improvement in the decision-making process by considering how much each new observation will reduce uncertainty about the objective function, not just the immediate improvement.

Another form of \emph{look-ahead} acquisition function is  \emph{Entropy Search}, which aims to reduce the uncertainty (quantified as the entropy) in the \emph{location} of the estimated optimum. 
The next point to be evaluated is the one that leads to the biggest decrease in the entropy of the estimated optimum location\shortcite{villemonteix2009informational,hennig2012entropy}:
    
   \begin{equation}
    \begin{aligned}
ES(\mathbf{x})=\mathrm{H}\bigl[\mathbf{x}^*|\mathcal{D}\bigl]- \mathbb{E}\Bigl[\mathrm{H}\bigl[\mathbf{x}^*|\mathcal{D}\cup\{(\mathbf{x},\mathbf{y}) \}\bigl]\Bigl]
    \end{aligned}
    \label{Background: equ:es}
    \end{equation}

\noindent In this expression,  $\mathrm{H}$ is the entropy function, $\mathbf{x}^*$ is the optimal solution of the problem, and $\mathcal{D}$ is the set of solutions observed so far. $ES$ provides a measure of the information gained on $x^*$  after observing the functions at point $\mathbf{x}$, resulting in the set of outputs $\mathbf{y}$.

   
Different authors have made attempts to extend EI to handle noisy settings \shortcite{ament2023unexpected,zhou2024corrected}. Most of these approaches, e.g., \shortciteN{letham2019cons}, rely on Monte Carlo methods, which sample from the posterior distribution to approximate the true expected improvement in the presence of noise. \shortciteN{quan2013simulation} proposed a modified version of the EI criterion, referred to as the \emph{modified expected improvement}:

\begin{equation}
\begin{aligned}
MEI(\mathbf{x})= \left[\mu_{sk}(x_{min})- \mu_{sk}(\mathbf{x})\right]\Phi \left(\frac{\mu_{sk}(x_{min})- \mu_{sk}(\mathbf{x})}{s(\mathbf{x})}\right)+ s(\mathbf{x}) \phi\left(\frac{\mu_{sk}(x_{min})- \mu_{sk}(\mathbf{x})}{s(\mathbf{x})}\right)
\end{aligned}
\label{equ:mei}
\end{equation}

In stochastic settings, this criterion improves the efficiency and effectiveness of the search compared with the traditional EI acquisition function. The MEI criterion is analogous to the $EI$ criterion, except for the use of $\mu_{sk}(x_{min})$ (the stochastic kriging predictor) instead of $f_{min}$ (the best-simulated value so far). This is a straightforward choice, as it provides an unbiased prediction given the heteroscedastic nature of the noise \shortcite{quan2013simulation,ankenman_nelson_staum_2010}. 

The \emph{Augmented Expected Improvement (AEI)}, proposed by \shortcite{huang2006global}, incorporates a noise-adjustment term:  

\begin{equation} \label{aei}
AEI(x) = \mathbb{E}\left[ \max\left( \mu(x_{min}) - \mu(x), 0 \right) \right] \times \left( 1 - \frac{\tau}{\sqrt{\tau^2 + s^2 (x)}} \right)
\end{equation}

\noindent where $\tau^2$ represents the variance of the noise. The goal is to avoid oversampling in very noisy regions. Note, however, that $AEI$ assumes homogeneous noise (some authors, e.g., \shortciteN{jalali2017comparison}, replace $\tau^2$ with $\tau^2(x)$ to reflect the presence of heteroscedastic noise). AEI performs best when the noise variance is known; in real-life applications, however, this is usually not the case. Researchers then often employ an additional Kriging model to estimate this variance \shortcite{wang2018constrained}.


Correlated Knowledge Gradient acquisition function in \shortciteN{frazier2009knowledge} is developed to handle heteroscedastic noise.
Empirical studies in \shortciteN{jalali2017comparison} and \shortciteN{picheny2013benchmark} show that Knowledge Gradient acquisition functions can be effective in noisy environments.

    


\vspace{2mm}
\noindent \textbf{Maximizing acquisition functions} is nontrivial, as they tend to be highly nonlinear and multi-modal \shortcite{Nguyen2023bayesian}. Heuristic algorithms such as adaptive random search \shortcite{zabinsky2003stochastic,zhigljavsky2012theory} or evolutionary methods \shortcite{tzanetos2020comprehensive,scott2011correlated} can be used, as well as local optimization approaches with multi-start (e.g., \shortciteN{loka2024cheap}). 
Multi-start methods enhance search efficiency by repeatedly restarting from new initial solutions, combining exploration of diverse regions with local search refinement. This approach helps avoid local optima and improves the chances of finding high-quality solutions. 
All of these are heuristic methods, which may require many function evaluations. Fortunately, this is usually not an issue as most acquisition functions are cheap to evaluate (exceptions are knowledge gradient and entropy search). The optimality of a heuristic approach cannot be guaranteed, though, and its performance typically depends on a number of user-defined parameters. To avoid these issues, some authors \shortcite{gonzalez2020multiobjective,jalali2017comparison,angun2023constrained} discretize the search space into a finite (but potentially large) space-filling set of so-called \emph{candidate points}. The acquisition function value is then evaluated for all unvisited alternatives, and the one with the highest value is chosen as the next infill point. This approach helps to ensure that differences in observed performance across algorithms are not caused by the heuristic at hand, but by differences in algorithmic approach. In high-dimensional problems, the use of candidate sets may become impractical, as the number of required points increases sharply with the number of dimensions.

While discretization is beneficial in studying and comparing different algorithms' performance, it has its downside. The optimization of a continuous problem by discretizing the solution space can sometimes lead to suboptimal solutions, regardless of the quality of the optimization algorithm. This arises from the fact that, no matter how extensive the discrete set is or how it is sampled, there is no assurance that this resulting set encompasses the global optimum of the problem. Consequently, it becomes impossible for any algorithm to discover this global optimum during the optimization process.

\section{Kriging vs. polynomial models: an illustrative example}


This section compares the performance of kriging and polynomial regression models in simulation optimization, using a simple $(R,s,S)$  inventory management system as an illustrative example. 

\textbf{Simulation Model.}  We adapted the $(R,s,S)$ inventory system described in \shortciteN{law2015simulation} to develop our simulation model. 
In this system, the inventory position ($I$) is monitored at the start of each month. When it falls to or below a specified reorder point ($s$), a replenishment order is placed such that the inventory position reaches the order-up-to level ($S$). Customer demands are met immediately if there are enough items in stock. If not, the excess demand is backlogged and is met as soon as there are enough units in inventory. The
time between demands is IID exponentially distributed with a mean of 0.1 months. The demand sizes are IID integer random variables; they follow a discrete distribution with values between 1 and 4 (occurring with probability $\frac{1}{6}$, $\frac{1}{3}$, $\frac{1}{3}$, and $\frac{1}{6}$ respectively).
Once a replenishment order has been placed, the time required for it to arrive (called the delivery lead time \shortcite{duran2023evaluating}) is random, following a uniform distribution between 0.5 and 1 month.

There is a fixed cost associated with placing an order (\$32).
Holding costs are incurred for each unit in stock at the end of each month (\$1 per item per month). Backlog costs arise when demand exceeds available stock, costing \$5 per unit per month. The system is evaluated over 120 months of operation. Further detail on the simulation model can be found in \shortciteN{law2015simulation}. The goal of the optimization is to determine the values $s$ and $S$ such that the expected total cost (consisting of holding, shortage, and order costs) is minimized. Each $(R,s,S)$ solution that needs to be evaluated is replicated 5 times in the simulation.

\textbf{Optimization algorithms:} Three different optimization algorithms are compared. The first algorithm, SK/MEI, is a kriging-based optimization algorithm using stochastic kriging (with a squared exponential kernel (Equation \ref{kernel se})) to fit the observed costs, and applying the MEI acquisition function (Equation \ref{equ:mei}) to select the next point to evaluate. The second algorithm, OK/EI, fits an ordinary kriging model (using the same kernel) to the observed sample means (thereby neglecting the noise on these sample means), and uses the EI acquisition function (Equation \ref{equ:ei}). 
The third algorithm, Polynomial Regression, fits a second-order polynomial to the sample means, with the next point chosen based solely on the predicted output. While SK/MEI and OK/EI use acquisition functions that account for metamodel uncertainty, Polynomial Regression selects points based on predicted values alone. 

We evaluate each algorithm on 25 macro-replications of the $(s,S)$ system, where each macroreplication starts from a different lhs design (consisting each time of 10 points).All three algorithms use the same initial design per macroreplication, to ensure a fair comparison. The candidate set of solutions consists of a two-dimensional grid of $(s,S)$ locations, with $s$ and $S$ integer values between 1 and 100; evidently, $s<S$ such that the grid consists of 4950 values in total. 
Each macroreplication ends after 100 infill points have been added. The observed solution with the smallest sample mean is then returned as the estimated optimal solution.

\begin{figure}[t!]
\centering
\includegraphics[width=\textwidth]{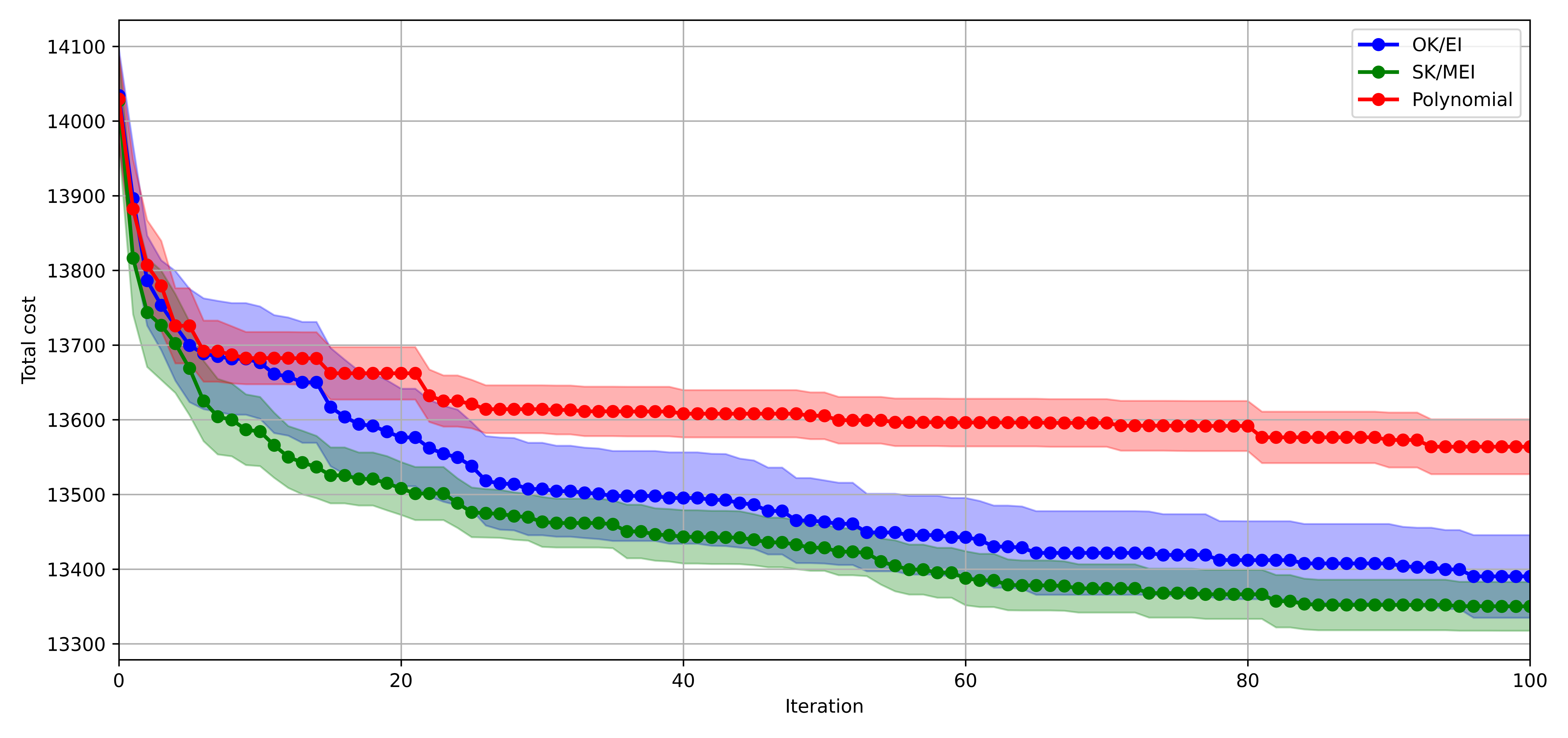}
\caption{Convergence comparison of three optimization Algorithms: SK/MEI, OK/EI, and Polynomial regression over 100 iterations, results of 25 macro-replications}
\label{fig: conv}
\end{figure}

\textbf{Discussion:} Figure \eqref{fig: conv} shows the evolution in the best expected total cost obtained, over the 100 sequential iterations, for each of the three algorithms. The lines show the average over the 25 macroreplications, and the shaded areas reflect the 95\% t-confidence intervals on this average.
SK/MEI, leveraging stochastic kriging and MEI to effectively balance exploration and exploitation, clearly outperforms OK/EI by focusing quickly on promising regions of the solution space.  This is not surprising, as OK/EI treats the noisy sample means as deterministic information, which may mislead the search. The solution quality obtained after 100 iterations by OK/EI is reached by SK/MEI within approx. $1/2$ of this budget. The noise handling within SK/MEI is key in improving data efficiency. 
Polynomial Regression shows the slowest convergence, and hence performs worst in terms of solution quality for this limited simulation budget.

\section{Extensions of kriging-based simulation optimization}

Thus far, we have focused exclusively on single-objective problems without constraints. However, many real-world optimization scenarios are more complex, involving multiple objectives simultaneously and/or constraints that limit the feasible solution space. These challenges have resulted in extensions of the basic  framework described above; in this section, we briefly explore these extensions.

\subsection{Multi-objective optimization}

In multi-objective optimization, the aim is to find solutions that simultaneously optimize multiple conflicting objectives. Without any preference information from the decision maker, the outcome is a set of Pareto-optimal solutions, where improving one objective means compromising at least one other objective \shortcite{hunter2019introduction}. A common method for approximating the Pareto frontier is by means of scalarization functions, which combine multiple objectives into a single value by assigning weights. The algorithm then attempts to obtain a high-quality approximation of the Pareto front by varying these weights.  The Augmented Tchebycheff function is often preferred over other scalarization methods, such as the Weighted Sum, because it better maintains diversity along the Pareto front and ensures that points on nonconvex regions of the front can be detected \shortcite{miettinen1999nonlinear}.

Multi-objective optimization in settings with noisy evaluations is challenging. Simply relying on noisy sample means to detect the front may lead to misclassification (truly Pareto-optimal solutions may be identified as dominated, while dominated solutions are considered Pareto-optimal). These errors prevent the correct identification of the truePareto-optimal solutions, especially when the replication budget is limited.

According to \shortciteN{rojas2020survey}, acquisition functions in kriging-based multi-objective simulation optimization fall into two categories. The first type consists of single-objective acquisition functions, originally designed for single-objective problems but adapted to multi-objective contexts by assessing improvement for each objective individually or through a scalarized version. The second type includes multi-objective acquisition functions, which measure the contribution of new points to the Pareto front using metrics such as hypervolume improvement, see, e.g., \shortciteN{yang2019multi}.

\subsection{Constrained optimization}

The key challenge in constrained kriging-based simulation optimization is ensuring that the algorithm accurately identifies feasible regions defined by (possibly stochastic) constraints. Failing to do so introduces a significant risk of returning infeasible or suboptimal solutions as the final outcome \shortcite{Aminiphdthesis}. In these problems, the selection of infill points is not only driven by how promising points are for improving the objective function, but also by how they may contribute to the estimation of the feasible region. 

Despite its importance, constrained kriging-based stochastic simulation optimization remains relatively underexplored \shortcite{amini2025constrained}. 
\shortciteN{wang2018constrained} developed an algorithm to address stochastic simulation optimization problems with only one constraint. Their acquisition function balances two factors: the point’s potential for improving the objective function and its likelihood of satisfying the constraint. This balance helps the algorithm in finding better solutions while learning the boundaries of the feasible region, reducing the risk of selecting infeasible solutions. \shortciteN{amini2023bayesian} extended this framework to handle problems with multiple constraints, incorporating stochastic kriging and a barrier function approach, which penalizes infeasible points and guides the search towards the optimal solution within the feasible region. Meanwhile, \shortciteN{angun2023constrained} integrates the Karush-Kuhn-Tucker (KKT) conditions into kriging-based simulation optimization. Their acquisition function multiplies the expected improvement by a cosine function, which measures how well a point satisfies the KKT conditions.

In constrained optimization, it is crucial to focus on enhancing feasibility estimations in areas where the objective function is estimated to have small values (in the case of minimization). This allows the algorithm to prioritize regions likely to contain the optimal solution, rather than spending resources in less relevant areas, making the algorithm more data-efficient. Both \shortciteN{wang2018constrained} and \shortciteN{amini2023bayesian} emphasize this approach, underscoring its importance for efficiency. In contrast, \shortciteN{angun2023constrained} take a different approach, aiming to enhance accuracy across the entire feasible region. While this improves overall precision, it requires a larger computational budget, making the algorithm less data-efficient.

\section{Conclusion}

In this tutorial, we explored the basics of kriging-based simulation optimization, emphasizing its data efficiency, particularly when handling stochastic problems with heteroscedastic noise. By comparing kriging-based methods to traditional polynomial approaches, we highlighted that kriging models are more flexible and better capture complex, non-linear relationships in I/O data. Furthermore, better uncertainty modeling in stochastic kriging leads to faster and more reliable convergence.

In industrial settings, data-efficient simulation optimization is essential for addressing complex, resource-intensive problems. Kriging-based methods provide a powerful framework for optimizing process or product designs within fewer simulation runs, reducing costs, and accelerating decision-making. This efficiency, coupled with the ability to handle noisy data, makes kriging-based optimization a vital tool for companies seeking to enhance competitiveness. As industries become increasingly reliant on simulation models and digital twins for product development and operational improvements, the use of data-efficient methods like kriging will continue to play a key role in driving innovation and sustaining growth.

\section*{ACKNOWLEDGMENTS}
The authors are grateful to the anonymous reviewers for their thoughtful feedback, which has enhanced the clarity of this work. This work was supported by the Flanders Artificial Intelligence Research Program (FAIR2) and the Research Foundation Flanders (FWO), grant number: G0A4624N.

\bibliographystyle{sw}
\bibliography{demobib}

\section*{AUTHOR BIOGRAPHIES}

\noindent {\bf Sasan Amini} is a postdoctoral researcher in the Computational Mathematics research group at Hasselt University, Belgium. He holds a Ph.D. in Business Economics, and his research primarily focuses on constraint handling methods in black-box optimization problems, with a particular emphasis on stochastic settings and problems involving expensive-to-evaluate functions.

\noindent {\bf Inneke Van Nieuwenhuyse} is a Professor at Hasselt University and KU Leuven, Belgium. Her research interests focus on (multi-objective) optimization of stochastic systems, especially in settings that are expensive to evaluate. To that end, she combines operations research techniques with machine learning approaches.

\end{document}